\documentclass{amsart}
\usepackage{amssymb,latexsym,enumerate}
\newtheorem{Thm}{Theorem}[section]

\newtheorem{Lem}[Thm]{Lemma}
\newtheorem{Cor}[Thm]{Corollary}
\newtheorem{Prob}[Thm]{Problem}
\newtheorem{Conj}[Thm]{Conjecture}
\numberwithin{equation}{section}

\DeclareMathOperator{\Q}{Q}
\DeclareMathOperator{\Loc}{L}

\DeclareMathOperator{\core}{core}

\DeclareMathOperator{\Stab}{Stab}
\DeclareMathOperator{\Homeop}{Homeo_+}

\DeclareMathOperator{\NS}{NS}
\DeclareMathOperator{\NF}{NF}

\DeclareMathOperator{\Aut}{Aut}
\DeclareMathOperator{\Fix}{Fix}
\DeclareMathOperator{\PLO}{PLO}
\DeclareMathOperator{\PLT}{PLT}
\DeclareMathOperator{\Cen}{\textbf{C}}
\DeclareMathOperator{\Zen}{\textbf{Z}}

\newcommand{\oc}[1]{(#1]}
\newcommand{\co}[1]{[#1)}
\begin{document}

\date{July 12, 2000}

\title[Left ordered groups]
{Left ordered groups with no nonabelian free subgroups}

\author[P. A. Linnell]{Peter A. Linnell}

\address{Math \\ VPI \\ Blacksburg \\ VA 24061--0123
\\ USA}

\email{linnell@math.vt.edu}
\urladdr{http://www.math.vt.edu/people/linnell/}

\begin{abstract}
There has been interest recently concerning when a left ordered group
is locally indicable.  Bergman and Tararin have shown that not all
left ordered groups are locally indicable, but all known examples
contain a nonabelian free subgroup.  We shall show for a large class
of groups not containing a nonabelian free subgroup, that any left
ordered group in this class is locally indicable.  We shall also show
that certain free products with an amalgamated cyclic subgroup are
left orderable.
\end{abstract}



\maketitle

\section{Introduction} \label{Sintroduction}

A group $G$ is left ordered if it has a total ordering $\le $ such
that $x \le y \Rightarrow gx \le gy$ whenever
$g,x,y \in G$.  For much
information on left ordered groups,
see the books \cite{kopytovmed,mura}.
Of course we say that a group $G$ is \emph{right} ordered
if it has a total ordering $\le $ such that $x \le y \Rightarrow xg
\le yg$ whenever $g,x,y \in G$. However using the involution $g
\mapsto g^{-1}$ of $G$, it is easy to see that a group is right
orderable if and only if it is left orderable.

Recall that a group is \emph{locally indicable} if and only if every
finitely generated subgroup $\ne 1$ has an infinite cyclic quotient.
Every locally indicable group is left orderable
\cite[theorem 7.3.1]{mura}, but the converse is not true, as has been
shown by Bergman \cite{bergman} and
Tararin \cite{tararin_ex}.  On the
other hand Chiswell and Kropholler \cite[theorem~A]{chiskrop} showed
that a solvable-by-finite left ordered group is locally indicable;
also Tararin \cite[theorem~3]{tararin_sol} has proved
that if $A\lhd G \ne 1$
are groups with $G/A$ finitely generated and solvable, $A$ abelian
and $G$ left orderable, then $G$ has a quotient isomorphic to an
infinite subgroup of $\mathbb Q$.   Further results in this direction
were obtained in \cite{rhem}.
In \cite{order} it was proved that
every left ordered elementary amenable group
is locally indicable, and
the question was raised of whether every left ordered amenable group
is locally indicable.
Let $\NF$ denote the class of groups
which contain no nonabelian free subgroup.
We shall consider the following stronger statement.
\begin{Conj} \label{Cmain}
A left ordered $\NF$-group is locally indicable.
\end{Conj}

We shall now describe the class of groups for which we shall prove
Conjecture \ref{Cmain}.
Let $\mathcal {P}$ denote the group of piecewise linear orientation
preserving self homeomorphisms of the unit interval $[0,1]$ with
multiplication defined as composition of functions.  Thus if $f,g \in
\mathcal {P}$, then $f$ is differentiable at all but a finite number
of points, and $(fg)(x) = f(g(x))$ for all $x \in [0,1]$.  Also let
$\NS$ denote the class of groups which have no nonabelian free
subsemigroup.  Now define $\mathcal {C}$ to be the smallest class of
groups which contains $\NS$ and $\mathcal {P}$, and is closed under
taking subgroups, homomorphic images,
group extensions and directed unions.  Clearly
$\mathcal{C}$ contains all elementary amenable groups and in
particular all solvable by finite groups, and it is not
difficult to show that $\mathcal {C} \subseteq \NF$
(see Corollary \ref{Cnf}).  Moreover $\mathcal {C}$ contains groups
which are not elementary amenable, such as the ubiquitous
Thompson's group (see \cite{floyd} for more information on this
topic).  Presumably not every $\NF$-group lies in the
class $\mathcal {C}$, though I know of no explicit example in the
literature.  We can now state
\begin{Thm} \label{Ttheorem}
A $\mathcal{C}$-group is left orderable if and only if it is
locally indicable.
\end{Thm}
Of course the result that $G$ is locally indicable (whether or not $G
\in \mathcal {C}$) implies that $G$ is left
orderable has already been
noted above.  For the reverse implication, we prove a stronger result
Theorem \ref{Tmain} which states that if $F \ne 1$ is a finitely
generated left orderable group and $F \supseteq G \in \mathcal{C}$,
then there exists a left-relatively convex subgroup $H \ne F$
(see Section \ref{Snotation})
such that $H \cap G \lhd G$, and $G/H \cap G$ has a self centralizing
torsion free normal abelian subgroup $A/H \cap G$ such that $G/A$ is
torsion free abelian.  Thus in the special case $F= G$ (so $G$ is
finitely generated and $\ne 1$), we see that $G$ has a quotient
isomorphic to $\mathbb {Z}$.

In Section \ref{Scircle} we shall use Theorem \ref{Ttheorem} to
prove the following result about
$\Homeop (S^1)$, the group of orientation preserving
homeomorphisms of the circle.
\begin{Cor} \label{Ccircle}
Let $G$ be a finitely generated subgroup of $\Homeop (S^1)$ such that
$G \in \mathcal {C}$.  Then
\begin{enumerate} [\normalfont (i)]
\item If $G$ is finite, then $G$ is cyclic.
\item If $G$ is infinite, then there exists $K \lhd H \lhd G$ such
that $G/H$ is cyclic and $H/K \cong \mathbb {Z}$.
\end{enumerate}
\end{Cor}
This should be compared with \cite[theorem 1.1]{shalen}, where
by considering the smaller group of orientation preserving
\emph{$C^{\infty}$ diffeomorphisms} of $S^1$, similar
but stronger results were obtained.

In Section \ref{Samalg} we shall use some of the techniques in this
paper to show that certain free products with amalgamation are
left orderable.  For example, we shall show in Theorem \ref{Tamalg}
that the free product of a left orderable group and a
torsion free nilpotent group with an amalgamated cyclic subgroup is
left orderable.
In the final section we shall briefly consider some examples of
left ordered groups which are not locally indicable.

Part of this work was carried out while I was at the
Sonderforschungsbereich in M\"unster.  I would like to thank
Wolfgang L\"uck for organizing my visit to M\"unster, and the
Sonderforschungsbereich for financial support.

\section{Notation, Terminology and Assumed Results}
\label{Snotation}

As usual $\mathbb {Q}$, $\mathbb {R}$, $\mathbb {Z}$ and $\mathbb
{N}$ will denote the rational numbers, real numbers, integers
and natural numbers $\{1,2, \dots \}$ respectively.  We
shall use the notation $G'$ for the commutator subgroup of the group
$G$, and if $g,x \in G$ and $X \subseteq G$, then $x^g = gxg^{-1}$,
$X^g = gXg^{-1}$, $\Cen_G(X) = \{g \in G \mid x^g = x$ for all $x
\in X\}$, $\Cen_G(x) = \Cen_G(\{x\})$,
and $\langle X \rangle$ denotes the subgroup generated by
$X$.  Also if $H \leqslant G$,
then $\core_G(H) = \bigcap_{g\in G}
H^g$, the largest normal subgroup of $G$ contained in $H$.

All mappings will be written on the left, in particular all group
actions will have the group acting on the left of the set.  If $G$ is
acting on a set $Y$ and $Z \subseteq Y$, then $\Stab_G(Z)$
will always
denote the pointwise stablizer of $Z$ in $G$: thus $\Stab_G(Z) =
\{g\in G \mid gz = z$ for all $z \in Z\}$, and we write $\Stab_G(y)$
for $\Stab_G(\{y\})$ when $y\in Y$.  Also if $H \subseteq G$, then
$\Fix_Y(H)$ is the fixed points of $H$, that is $\{y\in
Y \mid hy = y$ for all $h \in H\}$, and when $g \in G$ we write
$\Fix_Y(g)$ for $\Fix_Y(\{g\})$.  Then
obviously $\Fix_Y(X) = \Fix_Y(\langle X \rangle)$ whenever $X
\subseteq G$.

A totally ordered set $X$ is a set with a binary relation $\le$
such that for $x,y,z \in X$, either $x \le y$ or $y \le x$, $x\le y $
and $y \le x$ implies $x = y$, and $x \le y \le z$ implies $x \le z$.
Given totally ordered sets $X$ and $Y$, the map $\theta \colon X \to
Y$ is said to be order preserving if $x < y$ implies $\theta x <
\theta y$ whenever $x,y \in X$.  We shall let $\Aut(X)$ denote the
group of all order preserving permutations $X \to X$.  Note that if
$\theta$ is an order preserving bijection $X \to Y$, then $\theta
^{-1}$ is also order preserving and thus $\Aut(X)$ is indeed a group.
Also if $X \subseteq \mathbb {R}$ and $X$ is given the order induced
by the natural order on $\mathbb {R}$, then the elements of $\Aut(X)$
are homeomorphisms of $X$.

If $(G,\le)$ is a left ordered group and $K \leqslant G$, then we say
that $K$ is a convex subgroup of $G$ if $g \in G$, $j, k \in K$ and
$j \le g \le k$ implies $g \in K$.  In this case the left cosets of
$K$ in $G$, which we denote by $G/K$, is naturally
a totally ordered set
under the definition $gK < hK$ if and only if $g < h$
and $gK \ne hK$, for $g,h \in
G$.  Furthermore $G$ then acts as order preserving permutations on
$G/K$ according to the rule $g(hK) = ghK$.  We say that a subgroup
of $G$ is a \emph{left-relatively convex} subgroup
\cite[p.~127]{kopytovmed} if it is convex with respect to some left
order on $G$.  Conversely suppose $G$
acts faithfully as order preserving permutations on some totally
ordered set $X$.  Then, as described in
\cite[theorem 7.1.2]{mura}, we
can make $G$ into a left ordered group as follows.  Well order $X$,
and then for $f,g \in G$ with $f \ne g$, we say that $f < g$ if and
only if $f(x) < g(x)$ where $x$ is the least
element of $X$ such that
$f(x) \ne g(x)$.  Note that if $y \in X$ and $Y$ is the set of all
elements less than $y$, then $\Stab_G(Y)$ is a
convex subgroup of $G$
under this order and consequently $\Stab_G(Y)$
is a left-relatively
convex subgroup.  Therefore $\Stab_G(Y_0)$ is a
left-relatively convex
subgroup of $G$ for any subset $Y_0$ of $X$.
We need the following basic results about
left-relatively convex subgroups.
\begin{Lem} \label{Lrelconvex}
Let $G$ be a left ordered group, let $H$ be a
normal convex subgroup
of $G$, and let $\mathcal {B}$ be a set of left-relatively convex
subgroups of $G$.  Then
\begin{enumerate}[\normalfont (i)]
\item $\bigcap_{B \in \mathcal {B}} B$ is a left-relatively convex
subgroup of $G$.

\item If $\mathcal {B}$ is totally ordered by inclusion, then
$\bigcup_{B \in \mathcal {B}} B$ is a
left-relatively convex subgroup
of $G$.

\item If $B/H$ is a left-relatively convex
subgroup of $G/H$, then $B$
is a left-relatively convex subgroup of $G$.
\end{enumerate}
\end{Lem}
\begin{proof}
For (i) see \cite[proposition 5.1.10]{kopytovmed} or \cite[lemma
2.2(i)]{order}.  For (ii) see
\cite[proposition 5.1.7]{kopytovmed} or
\cite[lemma 2.2(ii)]{order}.  Finally for (iii), see \cite[lemma
2.1]{order}.
\end{proof}

We define $\PLO$ to be the class of groups which act
faithfully as piecewise linear orientation preserving self
homeomorphisms of $[0,1]$, and $\PLT$ the class
of groups which act
faithfully as piecewise linear orientation preserving self
homeomorphisms of $[0,1]$ which do not have a
common fixed point in $(0,1)$.
Thus $\PLT \subset \PLO$,
$G \in \PLO$ if and only if $G$ is isomorphic to a subgroup of
$\mathcal {P}$, and $G \in \PLT$ if and only
if $G$ acts faithfully as
piecewise linear orientation preserving self homeomorphisms of
$[0,1]$, and given $\epsilon > 0$ and $x \in (0,1)$,
there exist $f,g
\in G$ such that $f(x) < \epsilon$ and $g(x) > 1-\epsilon$.

Finally in this section, we need the following
refinement of the well
known fact that a countable left ordered group can be considered
as a subgroup of $\Aut(\mathbb {R})$ (see for example,
\cite[lemma 2.2]{witte}).
\begin{Lem} \label{Laction}
Let $G$ be a countable left ordered group, and let $H$ be a convex
subgroup of $G$ such that $H \ne G$.  Then there is an order
preserving action of $G$ on $\mathbb {R}$ with kernel $\core_G(H)$
such that $\Stab_G(0) = H$ and $\Stab_G(v) \ne G$ for all $v
\in \mathbb {R}$.
\end{Lem}
\begin{proof}
This follows from \cite[lemmas 2.4 and 2.3]{order}.
\end{proof}

\section{Extension Closed Classes of Groups}

Very similar results to the next lemma have been
proved before, see
for example \cite[lemma 3.1]{order}.  We shall prove a more general
result, which will hopefully avoid the need for further similar
results.  If  $\mathcal {D}$ is a class of groups which is closed
under taking subgroups, then we shall define
$\bar{\mathcal{D}}$ to be
the smallest class of groups containing $\mathcal{D}$ which is closed
under group extension and is closed under directed unions.
For arbitrary classes of groups $\mathcal {X}$ and $\mathcal{Y}$, we
shall let $H \in \Loc \mathcal {X}$ mean that every finite subset of
the group $H$ is contained in an $\mathcal{X}$-subgroup,
$H \in \Q \mathcal {X}$ to mean that $H$ is isomorphic to a quotient
group of an $\mathcal {X}$-group, and
$H \in \mathcal {X} \mathcal {Y}$ mean that $H$ has a normal
$\mathcal {X}$-subgroup $X$ such that $H/X \in \mathcal {Y}$.
If $\mathcal {X}$ is subgroup closed,
then $H \in \Loc \mathcal{X}$ if and
only if every finitely generated subgroup of $H$ is an $\mathcal
{X}$-group.
For each ordinal $\alpha$, the class of groups $\mathcal
{D}_{\alpha}$ is defined inductively by $\mathcal {D}_0 = \{1\}$,
$\mathcal {D}_{\alpha + 1} =
(\Loc \mathcal{D}_{\alpha}) \mathcal{D}$
and $\mathcal{D}_{\beta} = \bigcup_{\alpha < \beta} \mathcal
{D}_{\alpha}$ if $\beta$ is a
limit ordinal.  Setting $\mathcal {X} = \bigcup_{\alpha \ge 0}
\mathcal {D}_{\alpha}$, we can state

\begin{Lem} \label{Lelementary}
\begin{enumerate}[\normalfont (i)]
\item
Each $\mathcal {D}_{\alpha}$ and $\mathcal {X}$ is subgroup closed.

\item
If $\mathcal {D}$ is quotient group closed, then each $\mathcal
{D}_{\alpha}$ and $\mathcal {X}$ is also quotient group closed.

\item
$\mathcal {X} = \bar {\mathcal {D}}$.
\end{enumerate}
\end{Lem}
\begin{proof}
(i) This is easily proved by induction on $\alpha $, using the
fact that $\mathcal {D}$ is subgroup closed.
\smallskip

(ii) This is also easily proved by induction on $\alpha$, using the
fact that $\mathcal {D}$ is quotient group closed.

(iii) Clearly $\mathcal {X} \subseteq \bar
{\mathcal {D}}$, $\mathcal {X}
\supseteq \mathcal {D}$, and $\mathcal {X}$ is closed under directed
unions.  Therefore we need to prove that $\mathcal {X}$
is extension closed.

We show by induction on $\beta $ that $\mathcal {D}_{\alpha}
\mathcal {D}_{\beta} \subseteq \mathcal {D}_{\alpha +\beta }$;
the case $\beta = 0$ being obvious.  If
$\beta = \gamma +1$ for
some ordinal $\gamma $, then
\begin{align*}
\mathcal {D}_{\alpha} \mathcal {D}_{\beta}  &= \mathcal {D}_{\alpha}
((\Loc \mathcal {D}_\gamma ) \mathcal {D})
\subseteq (\mathcal {D}_{\alpha} (\Loc \mathcal {D}_\gamma ))
\mathcal {D} \subseteq (\Loc
(\mathcal {D}_{\alpha}  \mathcal {D}_\gamma )) \mathcal {D} \\
&\subseteq (\Loc \mathcal {D}_{\alpha +\gamma }) \mathcal {D}
\quad \text{(by induction)}  \\
&= \mathcal {D}_{\alpha +\beta }.  \\
\intertext{On the other hand if $\beta $ is a limit ordinal, then
$\mathcal {D}_{\beta}  = \bigcup_{\gamma  < \beta }
\mathcal {D}_\gamma $ and}
\mathcal {D}_{\alpha}  \mathcal {D}_{\beta}  &=
\mathcal {D}_{\alpha}  \left( \bigcup_{\gamma <\beta }
\mathcal {D}_\gamma
\right) = \bigcup_{\gamma <\beta }
\mathcal {D}_{\alpha}  \mathcal {D}_\gamma  \\
&\subseteq \bigcup_{\gamma <\beta }
\mathcal {D}_{\alpha +\gamma } \quad
\text{(by induction)} \\
&\subseteq \mathcal {D}_{\alpha +\beta }
\end{align*}
as required.
\end{proof}

For the rest of this paper, we let
$\mathcal {D} = \NS \cup \Q(\PLO)$.   Then clearly
$\bar {\mathcal{D}} \subseteq \mathcal {C}$ and $\mathcal{D}$ is
closed under taking subgroups, quotient groups, consequently
$\mathcal {C} = \bar {\mathcal {D}} = \mathcal{X}$.

\section{Proof of the Main Theorem} \label{Smain}

The statement and proof of the next lemma is just a reformulation of
\cite[assertion 2.1]{solodov_circ}.
\begin{Lem} \label{Lnf1}
Let $W$ be a nonempty closed subset of $\mathbb{R}$ and let $\alpha,
\beta, z \in \Aut(W)$.  Suppose $\Fix_W(z) = \emptyset$ and
$z$ commutes with $\alpha$ and $\beta$.  If $\Fix_W(\alpha) \ne
\emptyset \ne \Fix_W(\beta)$
and $\Fix_W(\alpha) \cap \Fix_W(\beta) =
\emptyset$, then $\langle \alpha, \beta \rangle $ contains a
nonabelian free subgroup.
\end{Lem}
\begin{proof}
Without loss of generality we may assume that
$0,1 \in W$, $\beta (0)
= 0$, and $z(0) = 1$.  Then $\beta (1) = 1$.  Suppose that
$\Fix_W(\alpha) \ne \emptyset$
yet $\Fix_W(\alpha) \cap \Fix_W(\beta) =
\emptyset$.  We may write $\mathbb
{R} \setminus \Fix_W(\alpha)$ and $\mathbb {R} \setminus
\Fix_W(\beta)$ as a disjoint union of open intervals, which we shall
call $\mathcal{A}$ and $\mathcal{B}$ respectively.
Then a finite
number $n$ of these intervals will cover $[0,1]$; let these intervals
be $(a_0, a_1)$, $(b_1, b_2)$, $(a_2,a_3)$, $(b_3, b_4)$, \dots,
$(b_{n-2}, b_{n-1})$, $(a_{n-1}, a_n)$ (so $n$ is an odd integer)
where the $(a_{2i}, a_{2i+1})$
are intervals in $\mathcal{A}$, the $(b_{2i+1}, b_{2i+2})$ are
intervals in $\mathcal{B}$, and
$a_0 < 0 < a_1  \le a_2 < a_3 \le \dots
\le a_{n-1} < 1 < a_n$, $0 < b_1 < b_2 \le b_3 < b_4
< \dots < b_{n-1} < 1$.  Note that $z a_0 = a_{n-1}$
and $z a_1 = a_n$,
and  $a_i \in \Fix_W(\alpha)$ and $b_i \in
\Fix_W(\beta)$ for all $i$.  Also $(b_i,a_i) \cap W \ne \emptyset$ if
$i$ is odd, and $(a_i, b_i) \cap W \ne \emptyset$ if $i$ is even.  To
see this let us consider the former case.  We have $b_i \in
\Fix_W(\beta)$, so certainly $b_i \in W$.  Also $b_i \notin
\Fix_W(\alpha)$, so by replacing $\alpha$ with $\alpha^{-1}$ if
necessary, we may assume that $\alpha (b_i) > b_i$.
Then $\alpha (b_i)
\in W$ and since $\alpha (a_i) = a_i$, it follows that
$\alpha(b_i) \in (b_i, a_i)$.  Similarly if $i$
is even, we can show that
$(a_i, b_i) \cap W \ne \emptyset$.  Now choose
$x_i \in (b_i,a_i) \cap
W$ if $i$ is odd ($1 \le i \le n-2$), and $x_i \in (a_i,b_i) \cap W$
($2 \le i \le n-1$) if $i$ is even.  Finally set $x_0 =
z^{-1}x_{n-1}$.

Set $P_1 = (x_0,x_1) \cup (x_2,x_3) \cup
\dots\cup  (x_{n-3},x_{n-2})$
and $Q_1 = (x_1,x_2) \cup (x_3,x_4)\cup \dots \cup
(x_{n-2},x_{n-1})$, and
for $r \in \mathbb {Z}$ define
\begin{align*}
z^rP_1 &= (z^r x_0,z^r x_1) \cup (z^r x_2,z^r x_3) \cup \dots\cup
(z^r x_{n-3},z^r x_{n-2})  \\
z^rQ_1 &= (z^r x_1,z^r x_2) \cup (z^r x_3,z^r x_4) \cup\dots \cup
(z^r x_{n-2},z^r x_{n-1}).
\end{align*}
Now set $P = \bigcup_{r \in \mathbb {Z}} z^rP_1$
and $Q = \bigcup_{r \in \mathbb {Z}} z^rQ_1$.
Observe that $P\cap Q = \emptyset$.  Indeed if $y \in P\cap Q$, then
by translating by $z^r$ for suitable $r$, we may assume that $y \in
(0,1)$, and then the result is clear.

If $i$ is even, then $(x_i, x_{i+1}) \subset (a_i,a_{i+1})$.
Since $\Fix_W(\alpha) \cap (a_i,a_{i+1}) = \emptyset$,
we see that either $\alpha (x) > x$
for all $x \in (a_i,a_{i+1}) \cap W$, or $\alpha (x) < x$
for all $x \in (a_i,a_{i+1}) \cap W$; without loss of generality we
may assume that $\alpha (x) > x$.
Now $a_i, a_{i+1} \in \Fix_W(\alpha)$,
$a_i < x_i < b_i$ and $b_{i+1} < x_{i+1} < a_{i+1}$,
hence there exists a positive integer $p_i$ such that $\alpha^r
(x_i, x_{i+1}) \subset (x_{i+1}, a_{i+1})$ and $\alpha^{-r} (x_i,
x_{i+1}) \subset (a_i, x_i)$ for all $r > p_i$.
Let $p$ be the maximum of the
of the $p_i$ ($0 \le i \le n-3$).  Since $(x_{i+1}, a_{i+1})$,
$(a_i, x_i) \subset Q$, we see that $\alpha^{pr} (x_i, x_{i+1})
\subset Q$ for all $i$ and for all $r \ne 0$, and it follows that
$\alpha^{pr} P\subset Q$ for all $r \ne 0$.  Similarly there exists a
positive integer $q$ such that
$\beta^{qr} Q \subset P$ for all $r \ne
0$.  It now follows from Klein's Table Tennis lemma
\cite[p.~130]{delaharpe} that $\langle
\alpha^p, \beta^q \rangle $ is a free group.
\end{proof}

\begin{Lem} \label{Lnf2}
Let $W$ be a nonempty closed subset of $\mathbb {R}$
and let $G,Z \leqslant
\Aut(W)$.  Suppose $G \in \NF$,
$Z$ centralizes $G$ and $\Fix_W(Z) = \emptyset$.
Let $H = \{g \in G \mid \Fix_W(g) \ne \emptyset \}$.
Then $G' \subseteq H \lhd G$ and
$\Fix_W(F) \ne \emptyset$ for
every finitely generated subgroup $F$ of $H$.
\end{Lem}
\begin{proof}
Let $\{f_1, \dots, f_n\}$ be a finite subset of $H$, and set $F =
\langle f_1, \dots, f_n \rangle$.  The result will
follow if we can prove that $\Fix_W(F) \ne \emptyset$, because then
clearly $H \lhd G$, and $G/H$ is abelian by \cite[Theorem
2.1]{beklaryan}.
Choose $a_i \in \Fix_W(f_i)$, and then select $a \in W$
such that $a < a_i$ for all $i$.
Now choose $z \in Z$ such that $za > a_i $ for all $i$.

Suppose the sequence $\{z^ra \mid r > 0 \}$ is bounded above,
and let $L$ be the least upper bound of the sequence.  Then $z^r a_i
\in \Fix_W(f_i)$ for all $i$ and $L$ is the least upper bound of the
sequence $\{z^r a_i \mid r > 0 \}$.  We conclude that $L \in
\Fix_W(f_i)$ for all $i$ and hence
$L \in \Fix_W(F)$.  Therefore we may
assume that the sequence $\{z^ra \mid r > 0\}$ is
not bounded above, and
similarly we may assume that the sequence
$\{z^ra \mid r < 0\}$ is not
bounded below.

Therefore we may assume that $\Fix_W(z) = \emptyset$.  It now follows
from Lemma \ref{Lnf1} that $\Fix_W(f) \ne \emptyset$
for all $f \in F$.
Since $z$ commutes with all elements of $F$, we see that every
element of $F$ fixes a point in $\co{a,za}$
and we deduce that the set $\{fa \mid
f \in F\}$ is bounded above by $za$.  If $M$ is the supremum of
the set $\{fa \mid f\in F\}$, then $M \in \Fix_W(F)$
and the result is
proven.
\end{proof}

The statement and proof of the next lemma is just a reformulation of
\cite[lemma 3.1]{solodov_line}
\begin{Lem} \label{Lns1}
Let $W$ be a nonempty closed subset of $\mathbb{R}$, let $\alpha,
\beta \in \Aut(W)$, let $a \in \Fix_W(\alpha)$,
and let $b \in \Fix_W(\beta)$.
Suppose that $a < b$ and
$\Fix_W(\alpha ) \cap \oc{a,b}
= \emptyset = \Fix_W(\beta) \cap \co{a,b}$.
Then $\langle \alpha,
\beta \rangle$ contains a nonabelian free subsemigroup.
\end{Lem}
\begin{proof}
By replacing $\alpha $ and or $\beta$ with their inverses if
necessary, we may assume that $\alpha (b) < b$ and $\beta (a) > a$.
Set $x = \beta (a)$ and note that $a < x < b$.  Since
$\beta (b) = b$, $\beta $ has no fixed points on $\co{a,b}$,
and $\beta (a) > a$, we see that
there exists a positive integer $n$ such that $\beta^n
\co{a,b} \subseteq (x,b)$.  Similarly
there exists a positive integer $m$ such that $\alpha^m
\oc{a,b} \subseteq (a,x)$.

We now show that the subsemigroup generated by $\alpha^m$ and
$\beta^n$ is free on those generators.  Set $\gamma = \alpha^m$ and
$\delta = \beta^n$.  Suppose to the contrary that two nontrivial
distinct finite products $\pi, \rho$
of the form $\ldots \gamma^{n_1} \delta^{n_2}
\gamma^{n_3} \delta ^{n_4} \ldots $, where the $n_i$ are positive
integers, yield the same element of $\Aut(W)$.  By cancelling on the
left, we may assume without loss of generality that $\pi = \gamma
\pi_1$ and $\rho = \delta \rho_1$ or 1,
where $\pi_1, \rho_1$ are also
products of the form $\ldots \gamma^{n_1} \delta^{n_2}
\gamma^{n_3} \delta ^{n_4} \ldots $.  Since $\pi_1 x , \rho_1 x
\in (a,b)$, we see that $\gamma \pi_1 x \in (a,x)$ and
$\delta \rho_1 x$ or $1x \in \co{x,b}$.
Thus $\pi x \ne \rho x$ and we have a
contradiction.  We deduce that the subsemigroup generated by
$\alpha^m$ and $\beta^n$ is free on those generators and the
result follows.
\end{proof}

The statement and proof of the next lemma is just a reformulation of
\cite[lemma 3.6]{solodov_line}
\begin{Lem} \label{Lns2}
Let $W$ be a nonempty closed subset of $\mathbb {R}$, let $n$ be a
positive integer, let $\alpha_1, \dots,
\alpha_n \in \Aut(W)$, and let
$G = \langle \alpha_1, \dots, \alpha_n \rangle$.
Suppose $G \in \NS$.
If $\Fix_W(\alpha_i) \ne \emptyset $ for all $i$, then $\Fix_W(G) \ne
\emptyset$.
\end{Lem}
\begin{proof}
For each $i$, we may write
$\mathbb {R} \setminus \Fix_W(\alpha_i)$ as
a disjoint union of open intervals,
say $\bigcup_j I_{ij}$, where each
$I_{ij}$ is an open interval.  Then $I_{ij} \ne \mathbb {R}$ for all
$i,j$, because $\Fix_W(\alpha_i) \ne \emptyset$.
Suppose $\Fix_W(G) = \emptyset$.  If $I_{ij} \subseteq I_{kl}$
and $(i,j) \ne (k,l)$ then
$i\ne k$, so we may choose $i,j$ such that $I_{ij}$ is not contained
in any other open interval.  Using the fact that $\mathbb {R}$ is
connected, we may now choose $k,l$ so that $I_{ij}$ has nonempty
intersection with $I_{kl}$, and also does not contain $I_{kl}$.
Clearly $i \ne k$.  Write $I_{ij} \cap
I_{kl} = (a,b)$, and assume without loss of generality that $a \in
I_{kl}$ and $b \in I_{ij}$.  The result now follows by applying
Lemma \ref{Lns1} with $\alpha = \alpha_i$ and $\beta = \alpha_k$.
\end{proof}

\begin{Lem} \label{Lplt}
Let $G \in \PLT$, and let
$A$ and $B$ be finitely generated subgroups of $G'$.
Then there exists $g \in G$ such that $A^g$ and $B$ centralize each
other.
\end{Lem}
\begin{proof}
We may view
$G$ as a subgroup of the piecewise linear orientation
preserving homeomorphisms of $[0,1]$ such that $G$ fixes no point in
$(0,1)$.
Define $H = \{g \in G \mid \text{there exists } \epsilon > 0
\text{ such that } g(t) = t \text{ for all } t \in  [0,\epsilon] \cup
[1-\epsilon, 1]\}$.  Then obviously $H \lhd G$ and we see that
$H \supseteq G'$.  Therefore if $A$ and $B$ are finitely
generated subgroups of $G'$, there exists $0 < r < s < 1$ such that
$c$ is the identity map outside $(r,s)$ for all $c \in A \cup B$.
Since $G$ does not fix any point in $(0,1)$, there exists $g \in G$
such that $gr > s$.  Then $g A g^{-1}$ fixes all points outside
$(gr,gs)$ and the result follows.
\end{proof}

\begin{Lem} \label{Lorbits}
Let $G \in \PLO$ be a finitely generated group.
Then there exists a series $G' = G_0 \supseteq G_1 \supseteq \dots
\supseteq G_n =1$ with the property that
$G_i \lhd G$ and if $A,B$ are
finitely generated subgroups of $G_{i-1}/G_i$, then there exists $g
\in G/G_i$ such that $A^g$ and $B$ centralize each other, for $i = 1,
\dots, n$.
\end{Lem}
\begin{proof}
Write $G = \langle g_1, \dots , g_m \rangle$ and consider $G$ as a
group of orientation preserving homeomorphisms of $[0,1]$.  Let $W'$
denote the complement $[0,1] \setminus W$ of a subset $W$ of $[0,1]$.
Since $\Fix_{[0,1]} (g_i)'$ is a finite union of open intervals and
$\Fix_{[0,1]}(G)' = \bigcup_i \Fix_{[0,1]} (g_i)'$, we see that
$\Fix_{[0,1]} (G)'$ is a finite union of open intervals.  Therefore
$\Fix_{[0,1]} (G)'$ is a finite union of disjoint open intervals, say
\[
(a_1,b_1) \cup (a_2,b_2) \cup \dots \cup (a_n,b_n)
\]
where $0 \le a_1 < b_1 \le a_2 < b_2 \le
a_3 < \dots \le a_n < b_n \le 1$.
Set $H_i = \Stab_G([a_i,b_i])$, $G_0 = G'$,
and $G_i = G_0 \cap H_1 \cap \dots \cap H_i$ for
$i = 1, \dots, n$.  Then $G' =G_0 \supseteq G_1 \supseteq \dots
\supseteq G_n = 1$ and $G_i \lhd G$ for all $i$.
Since $g[a_i,b_i] = [a_i,b_i]$ for all $g \in G$,
we see that $G/H_i \leqslant \Aut([a_i,b_i])$.

Now suppose $A/G_i$ and $B/G_i$ are finitely generated
subgroups of  $G_{i-1}/G_i$.  Then $AH_i/H_i$ and
$BH_i/H_i$ are finitely generated subgroups of $G'H_i/H_i$, so
by Lemma \ref{Lplt} there exists $g \in G$ such that $[A^g,B]$ (the
commutator of $A^g$ and $B$) is contained in $H_i$.  Therefore
$[A^g,B] \subseteq H_i \cap G_{i-1} = G_i$ and the result follows.
\end{proof}

\begin{Lem} \label{Lplo}
Let $1 \ne G \in Q(\PLO)$ and suppose $G$ is finitely generated.
Then there exists $1 < H \lhd G$ such that if
$A,B$ are finitely generated subgroups of $H$, then there exists $g
\in G$ such that $A^g$ and $B$ centralize each other.
\end{Lem}
\begin{proof}
Write $G = P/K$ where $P \in \PLO$ and $K \lhd P$.  If $G = \langle
Kp_1, \dots, Kp_d \rangle$, then replacing $P$ with $\langle p_1,
\dots, p_d \rangle $ and $K$ with $K \cap \langle p_1, \dots, p_d
\rangle $, we may assume that $P$ is finitely generated.  By Lemma
\ref{Lorbits}, there exists a series $P'=P_0 \supseteq P_1 \supseteq
\dots \supseteq P_n = 1$ such that $P_i \lhd P$ for all $i$ with the
property that if $C/P_i$ and $D/P_i$ are finitely generated subgroups
of $P_{i-1}/P_i$, then there exists $p \in P$ such that $C^p/P_i$ and
$D/P_i$ centralize each other.  If $P' \subseteq K$, then we may take
$H = G$ and we are finished.  Otherwise we may let
$m$ be the smallest integer such
that $P_m \subseteq K$ and set $H = P_{m-1}K/K$, a nontrivial normal
subgroup of $G$.  If $A$ and $B$ are finitely generated subgroups of
$H$, then there exist finitely generated subgroups $C$ and $D$ of
$P_{m-1}$ such that $A = CK/K$ and $B = DK/K$.  Then we can find $p
\in P$ such that $C^pP_m/P_m$ and $DP_m/P_m$ centralize each other.
If $g = pK$, then $A^g$ and $B$ centralize each other as required.
\end{proof}

\begin{Cor} \label{Cnf}
$\mathcal {C} \subseteq \NF$.
\end{Cor}
\begin{proof}
It is easy to see that $\NF$ is closed under taking subgroups,
quotient groups, group extensions and directed unions.  Since $\NS
\subseteq \NF$, it will now be sufficient to show that $\mathcal {P}
\in \NF$.
However if $G$ is the free group on two generators, $1 \ne
H \lhd G$, $x \in H \setminus H'$ and $y \in H'\setminus 1$,
then there is no $g \in
G$ such that $x^g$ and $y$ centralize each other.
We deduce from Lemma \ref{Lplo} that $G \notin \PLO$ and the
proof is complete.
\end{proof}

\begin{Lem} \label{Lpl}
Let $W$ be a nonempty closed subset of $\mathbb {R}$ and let $G
\leqslant \Aut(W)$.
Suppose $1 \ne G \in \Q(\PLO)$ and that $G$ is finitely generated.
Then there exists $H \lhd G$ such that $H \ne 1$
and $\Fix_W(F) \ne \emptyset$
whenever $F$ is a finitely generated subgroup of $H'$.
\end{Lem}
\begin{proof}
Using Lemma \ref{Lplo}, we can find $H \lhd G$ such that $H \ne 1$
with the property that if $A,B$ are finitely generated subgroups of
$H$, then there exists $g \in G$ such that $A^g$ and $B$ centralize
each other.  The result is obvious if $\Fix_W(A) \ne \emptyset$ for
all finitely generated subgroups $A$ of $H$, so we may assume that
there is a finitely generated subgroup $A$ of $H$ such that
$\Fix_W(A) = \emptyset$.  Let $B$ be any finitely generated subgroup
of $H$.  Then there exists $g\in G$ such that $A^g$ and $B$
centralize each other.  Since $B \in \NF$ by Corollary
\ref{Cnf}, we deduce from Lemma \ref{Lnf2} that $\Fix_W(E) \ne
\emptyset$ for every finitely generated subgroup $E$ of $B'$.  We
conclude that $\Fix_W(F) \ne \emptyset$ for every finitely generated
subgroup $F$ of $H'$ as required.
\end{proof}

\begin{Lem} \label{Lmetabelian}
Let $W$ be a nonempty closed subset of $\mathbb {R}$ and
let $H \lhd G \leqslant \Aut(W)$ with $H$ solvable.
Suppose either $G/H \in \NS$, or $G/H \in
\Q(\PLO)$ and $G/H$ is finitely generated.  Then
either $G'' = 1$, or there exists $F \lhd G$ such that $F\ne 1$ and
$\Fix_W(E) \ne \emptyset$ whenever $E$ is a finitely generated
subgroup of $F$.  Furthermore in the former case there exists $A \lhd
G$ such that $A$ and $G/A$ are torsion free abelian, and $\Cen_G(A) =
A$.
\end{Lem}
\begin{proof}
First suppose $G$ has a nontrivial normal abelian subgroup $A$.
Here we set $F = \{a \in A \mid \Fix_W(a) \ne \emptyset \}$.
Since $A \in \NS$, we see from Lemma \ref{Lns2} that
$F \leqslant A$ and $\Fix_W(E) \ne \emptyset$ for all finitely
generated subgroups $E$ of $F$, so we may assume that $F=1$.
Let $C = \Cen_G(A)$ and let
$B = \{b \in C \mid \Fix_W(b) \ne \emptyset \}$.
Since $G \in \NF$ by Corollary \ref{Cnf},
Lemma \ref{Lnf2} shows that $B$ is a normal
subgroup of $G$ and that $\Fix_W(E) \ne \emptyset $ for all finitely
generated subgroups $E$ of $B$.  Therefore we
may assume that $B = 1$.
Using \cite[lemma 4.1]{order}, we see that $C$ and $G/C$ are
abelian and it follows that $G'' = 1$.  Thus we may assume that $G$
has no nontrivial normal abelian subgroup, so in particular $H = 1$.

We now have two cases to consider, namely $G \in \NS$ and $G \in
\Q(\PLO)$
and is finitely generated.  In the former case the result
follows from Lemma \ref{Lns2} and \cite[lemma 4.1]{order},
while in the latter case the result
follows from Lemma \ref{Lpl}.
\end{proof}

\begin{Lem} \label{Lmain}
Let $H \lhd G \leqslant F$ be groups such that $F$ is finitely
generated.  Assume that $G/H$ has a solvable normal subgroup $K/H$
such that $G/K \in \mathcal {D}$ and is finitely generated.
Suppose $F$ is left
orderable and there exists a left-relatively
convex subgroup $B$ of $F$
such that $H \subseteq B \ne F$.  Then there exists a left-relatively
convex subgroup $B_1$ of $F$ such that $B_1 \ne F$, $B_1 \cap G \lhd
G$, and $G/B_1 \cap G$ has a self centralizing torsion free abelian
normal subgroup $B_2/B_1 \cap G$ such that $G/B_2$ is torsion free
abelian.
\end{Lem}
\begin{proof}
For each $X \subseteq F$, let $\mathfrak {c} X$ denote the smallest
left-relatively convex subgroup of $F$
containing $X$, and let $\mathcal
{S} = \{I \lhd G \mid H \subseteq I$ and $\mathfrak {c} I \ne F\}$.
Then $\mathcal {S}$ is partially ordered by inclusion.  Suppose
$\mathcal {T}$ is a nonempty chain in $\mathcal {S}$.
Then $\bigcup_{I \in \mathcal {T}} \mathfrak{c} I$
is a left-relatively
convex subgroup of $F$ by Lemma \ref{Lrelconvex}, which is not the
whole of $F$ because $F$ is finitely generated and $\mathfrak {c} I
\ne F$ for all $I \in \mathcal {T}$, consequently $\mathcal {T}$ is
bounded above by $\bigcup_{I \in \mathcal {T}} I$.  But $H \in
\mathcal {S}$ because $\mathfrak {c} H \subseteq B \ne F$,
hence $\mathcal {S}
\ne \emptyset$ and we may apply Zorn's lemma to deduce that $\mathcal
{S}$ has a maximal element $E$ say.  Set $B_1 = \bigcap_{g \in
G}(\mathfrak{c} E)^g$, which by Lemma \ref{Lrelconvex}
is a left-relatively
convex subgroup of $F$, so using the maximality of $E$ we see that
$B_1 \cap G = E$ (thus $B_1 = \mathfrak {c}E$).  If $G/E$ has a self
centralizing torsion free abelian normal subgroup $B_2/E$ such that
$G/B_2$ is torsion free abelian, then we are finished so we assume
that this is not the case.

By Lemma \ref{Laction}, there is an order preserving action of $F$ on
$\mathbb {R}$ with kernel $\core_F(\mathfrak {c} E)$
such that $\Stab_F(0)
= \mathfrak {c}E$, and $\Stab_F(v) \ne F$ for
all $v \in \mathbb {R}$.
Replacing $F$ with $F/\core_F(\mathfrak {c}E)$ and using Lemma
\ref{Lrelconvex}, we may assume that $\core_F (\mathfrak {c} E) = 1$.

Let $W = \Fix_{\mathbb {R}} (E)$.  Then $W$ is a nonempty closed
subset of $\mathbb {R}$, and $G/E$ is naturally a subgroup of $\Aut
(W)$.  Using the hypotheses of the Lemma, there is a normal solvable
subgroup $K_1/E$ of $G/E$ such that $G/K_1 \in \mathcal {D}$ and is
finitely generated.  By Lemma
\ref{Lmetabelian}, there is a nontrivial normal subgroup $A/E$ of
$G/E$ such that $\Fix_W(C) \ne \emptyset$ whenever $C$ is a finitely
generated subgroup of $A/E$.  Write $A = \bigcup_{i \in \mathbb {N}}
A_i$ where $E \leqslant A_1 \leqslant A_2 \leqslant
\cdots $ and $A_i/E$ is finitely generated
for all $i$ (if $A/E$ is finitely generated, we may choose $A_i = A$
for all $i$), and $X_i = \Fix_{\mathbb {R}}(A_i)$.  Then $X_i \ne
\emptyset$ for all $i$, and $\Stab_F (X_1) \leqslant \Stab_F (X_2)
\leqslant \cdots $ is an ascending chain of
left-relatively convex subgroups
of $F$ with the property that $A_i \subseteq \Stab_F(X_i) \ne F$ for
all $i$.  Furthermore $\bigcup_{i \in \mathbb {N}} \Stab_F(X_i)$
is a left-relatively convex
subgroup by Lemma \ref{Lrelconvex}, which cannot
be $F$ itself because
$F$ is finitely generated.  We deduce that $\mathfrak {c} A \ne F$
which contradicts the maximality of $E$ and finishes the proof.
\end{proof}

\begin{Thm} \label{Tmain}
Let $G \leqslant F \ne 1$ be groups such that
$G \in \mathcal {C}$ and
$F$ is finitely generated and left orderable.  Then there exists a
left-relatively convex subgroup $B$ of $F$ such that $B \ne F$, $B
\cap G \lhd G$, and $G/B \cap G$ has a self centralizing torsion free
abelian normal subgroup $A/B \cap G$ such that $G/A$ is torsion free
abelian.
\end{Thm}
\begin{proof}
We shall prove the result by transfinite
induction on $G$, so by Lemma
\ref{Lelementary} choose the least ordinal $\alpha$ such that $G \in
\mathcal {D}_{\alpha}$ and assume that the result is true whenever $H
\in \mathcal {D}_{\beta}$ and $\beta < \alpha$.
Now $\alpha$ cannot be
a limit ordinal, and the result is clearly true if $\alpha = 0$.
Therefore we may assume that $\alpha = \gamma + 1$ for some ordinal
$\gamma$, and then there exists $H \lhd G$ such
that $G/H \in \mathcal
{D}$ and $H \in \Loc \mathcal {D}_{\gamma}$.  Using Lemma
\ref{Lelementary}, we may write $H = \bigcup_{i \in \mathbb {N}} H_i$
where $H_1 \leqslant H_2 \leqslant \dots \leqslant H$ and every
subgroup of $H_i$ is in $\mathcal {D}_{\gamma}$
for all $i$.  For each $X
\subseteq F$, let $\mathfrak {c} X$ denote the
smallest left-relatively
convex subgroup of $F$ containing $X$.

First consider the case $G/H$ is finitely generated.
We have an ascending chain of left-relatively convex subgroups
$\mathfrak {c} (H_1'') \leqslant \mathfrak {c} (H_2'') \leqslant
\cdots$, so their union is also a left-relatively convex subgroup by
Lemma \ref{Lrelconvex} which contains $H''$.
The inductive hypothesis
shows that $\mathfrak {c} (H_i'') \ne F$ for all $i$ and since $F$ is
finitely generated, we deduce that $\mathfrak {c} (H'') \ne F$.  But
$G/H''$ has the solvable normal subgroup $H/H''$ such that $G/H \in
\mathcal {D}$, so the result follows from
an application of Lemma \ref{Lmain}.

Finally we need to consider the case $G/H$ is not finitely generated.
Here we write $G = \bigcup_{i \in \mathbb {N}} G_i$, where $H
\leqslant G_i \leqslant G$ and $G_i/H$ is finitely generated for all
$i$.  We now have an ascending chain of left
relatively convex subgroups
$\mathfrak {c} (G_1'') \leqslant \mathfrak {c} (G_2'') \leqslant
\cdots$, so their union is also a left-relatively convex subgroup by
Lemma \ref{Lrelconvex} which contains $G''$.  By the case $G/H$ is
finitely generated considered in the previous paragraph, we know
that $\mathfrak {c} (G_i'') \ne F$ for all $i$ and since $F$ is
finitely generated, we deduce that $\mathfrak {c} (G'') \ne F$.
Another application of Lemma \ref{Lmain} completes the proof.
\end{proof}

\section{Groups of homeomorphisms of the circle} \label{Scircle}

Proof of Corollary \ref{Ccircle}.
\begin{proof}
By \cite[lemma 2.3]{witte},
we may lift the action of $G$ on $S^1$ to an action of a group
$H$ on $\mathbb {R}$; specifically $H$ is a left orderable group
with a central subgroup $Z$ such that $Z \cong \mathbb {Z}$ and
$H/Z \cong G$.  Note that $H$ is finitely generated because $G$
is finitely generated.  If $G$ is finite, then $H$ is a torsion
free group with an infinite central cyclic subgroup
of finite index and it follows
that $H \cong \mathbb {Z}$.  We deduce that $G$ is cyclic.  Therefore
we may assume that $G$ is infinite.

By Theorem \ref{Ttheorem} $H$ has a normal subgroup $K$ such that
$H/K \cong \mathbb {Z}$.  If $K \cap Z \ne 1$, then $Z/K \cap Z$
is finite, consequently $KZ/K$ is a finite subgroup of $H/K$ and
we deduce that $H/KZ \cong \mathbb {Z}$.  It follows that
$G$ has an infinite cyclic quotient, so we may assume that $K
\cap Z = 1$.

Note that $H/KZ$ is a finite cyclic group.
Since $KZ$ has finite index in the finitely
generated group $H$, we see that $K$ is finitely generated.  Moreover
$K$ is infinite, so by Theorem \ref{Ttheorem} there exists $L \lhd K$
such that $K/L \cong \mathbb {Z}$.  But
\[
\frac{KZ}{LZ} \cong \frac{K}{(LZ) \cap K}
= \frac{K}{L(Z \cap K)} = \frac{K}{L}
\]
and the result follows.
\end{proof}

\section{Free products with amalgamation}
\label{Samalg}

\begin{Lem} \label{Lreverse}
Let $G \leqslant \Aut(\mathbb {R})$.  Then there is an action
$\alpha$ of $G$ on $\mathbb {R}$ by orientation preserving
homeomorphisms with the following properties.
\begin{enumerate}[\normalfont (i)]
\item
If $c \in G$ and
$c(r) > r$ for all $r \in \mathbb {R}$, then
$(\alpha c)r < r$ for all $r \in \mathbb {R}$.
\item
If $c \in G$ and
$c(r) < r$ for all $r \in \mathbb {R}$, then
$(\alpha c)r > r$ for all $r \in \mathbb {R}$.
\end{enumerate}
\end{Lem}
\begin{proof}
Define an action $\alpha$ of $G$ on $\mathbb {R}$ by
$(\alpha g) r = -g(-r)$ for $g \in G$.
This action has the required properties.
\end{proof}

\begin{Lem} \label{Lfixempty}
Let $G \leqslant \Aut(\mathbb {R})$, let $H$ be a left ordered
group, let $1 \ne c \in G$, let $C = \langle c \rangle$, and
let $1 \ne h\in H$.
Identify $C$ with $\langle h \rangle$ via the isomorphism
$c^n \mapsto h^n \colon C \to \langle h \rangle$ for $n \in
\mathbb {Z}$.
Suppose $\Fix_{\mathbb {R}}(c) = \emptyset$.  Then $G *_C H$ is
left orderable.
\end{Lem}
\begin{proof}
Write $H = \bigcup_i H_i$, where the $H_i$ are finitely generated
subgroups containing $h$.  Then $G *_C H = \bigcup_i G*_C H_i$, and
if each of the $G *_C H_i$ is left orderable, then so is $G *_C H$
by \cite[7.3.2]{mura}.  Therefore we may assume that $H$ is
finitely generated.  Using Lemma \ref{Laction}, we can view
$H$ as a subgroup of $\Aut (\mathbb {R})$.

We may write $\mathbb {R} \setminus \Fix_{\mathbb {R}} (h)$
as a countable disjoint union of nonempty open sets, say
$\bigcup_i P_i$.  On each $P_i$, either $h(x) > x$ for all
$x \in P_i$, or $h(x) < x$ for all $x \in P_i$.  Using Lemma
\ref{Lreverse}, for each $i$ there is an action of $G$ on $P_i$
by orientation preserving homeomorphisms with the property
that either $h(x) > x$ and $c(x) > x$ for all $x$, or
$h(x) < x$ and $c(x) < x$ for all $x$.
Then by \cite[theorem 10]{fine}
we may assume that $h = c$ on $P_i$.  We have
now defined an action of $G$ on $\bigcup_i P_i$, and we extend
this to an action
$\alpha$ on the whole of $\mathbb {R}$ by defining $\alpha
g$ to be the
identity on $\Fix_{\mathbb {R}}(h)$ for all $g \in G$.  Clearly
$\alpha (G) \subseteq \Aut (\mathbb {R})$ and $\alpha (G) \cong G$.
Thus we can define a group homomorphism $\theta \colon G *_C H
\to \Aut(\mathbb {R})$
by $\theta g = \alpha g$ for $g \in G$ and $\theta h = h$ for $h \in
H$, because $\alpha (c^n) = h^n$ for $n \in \mathbb {Z}$.
The result now follows from \cite[theorem 6.2.3]{kopytovmed}.
\end{proof}

\begin{Thm} \label{Tamalg}
Let $G$ be a left ordered group, let $H$ be a torsion free
nilpotent group, and let $C$ be a cyclic group.  Then $G *_C H$ is
left orderable.
\end{Thm}
\begin{proof}
If $C = 1$ then the result follows from
\cite[\S 2.4 on p.~37 and theorem 7.3.2]{mura}, so
we may assume that $C$ is infinite cyclic.  We will assume that $C$
is a subgroup of $H$ and write $C = \langle c \rangle$, where $1 \ne
c \in H$.  Let $1 \ne g \in G$ and identify $C$ with $\langle g
\rangle$ via the isomorphism $c^n \mapsto g^n$ for $n \in \mathbb
{Z}$.  We need to prove that $G *_C H$ is left orderable.

Write $H = \bigcup_i H_i$, where the $H_i$ are finitely generated
subgroups containing $C$.  Then $G *_C H = \bigcup_i G *_C H_i$, and
if each of the $G *_C H_i$ is left orderable, then so is $G *_C H$
by \cite[7.3.2]{mura}.  Therefore we may assume that $H$ is
finitely generated.  We shall use induction on the Hirsch length of
$H$ (so if $1 = H_0 \lhd H_1 \lhd \dots \lhd H_n = H$ is a normal
series for $H$ with $H_i/H_{i-1}$ infinite cyclic for all $i$, then
$n$ is the Hirsch length of $H$).

First suppose the Hirsch length of $H$ is 1.  This means that $H$ is
infinite cyclic, say $H = \langle h \rangle$ where $h$ has infinite
order.  Then we can view $H$ as a subgroup of $\Aut (\mathbb {R})$
by letting $H$ act on $\mathbb {R}$ according to the rule
$h (r) = r+1$ for all $r \in \mathbb {R}$.  Then $\Fix_{\mathbb {R}}
(c) = \emptyset$ and the result follows from Lemma \ref{Lfixempty}.
Therefore we may assume that the Hirsch length of $H$ is at least 2.

Let $Z$ be a nontrivial cyclic central subgroup of $H$ such that
$H/Z$ is torsion free.  Then $H/Z$ is left orderable because
$H/Z$ is a torsion free nilpotent group
\cite[\S 2.4 on p.~37]{mura}.

Suppose $c \notin Z$.  We have an epimorphism
$G *_C H \twoheadrightarrow G *_{CZ/Z} H/Z$.
Let $K$ be the kernel of this map.
Then $K \cap G = 1$ and $K \cap H = Z$,
so applying \cite[I.7.7]{dicksdun}
we see that $K$ is free
and consequently left orderable.  By induction $G *_{CZ/Z} H/Z$
is left orderable, and so the result follows from
\cite[7.3.2]{mura}.

Finally we need to consider the case $c \in Z$.
We have an epimorphism $G *_C H \twoheadrightarrow H/Z$.
Let $K$ be the kernel of this map, and let
$F = G *_C H$.  With $G *_C H$ we have
an associated standard tree $T$
\cite[I.3.4 definitions]{dicksdun},
and $F$ acts on this tree.  The vertices of $T$ are the
left cosets $fG$ and $fH$, and the edges are the
left cosets $fC$, where $f \in F$.  A fundamental $F$-transversal
\cite[2.6 proposition]{dicksdun} for $T$ consists of the vertices
$G,H$ and the edge $C$.  Let $X$ be a transversal for
$Z$ in $H$.  Then a fundamental
$K$-transversal $T_0$ for $T$ consists
of the vertices $xG$ and $xH = H$, and the edges
$xC$, where $x \in X$.
The stabilizers of the vertices of $T_0$
are of the form $G^x$ and $Z^x = Z$, and the stabilizers
of the edges are of the form $C^x = C$ for $x \in X$.
It follows that $K$ is the fundamental
group of a graph of groups \cite[I.3.4 definitions]{dicksdun}
of the following form

\vspace{-3mm} %This \vspace is to eliminate a bad page break
{\large
\setlength{\unitlength}{1mm}
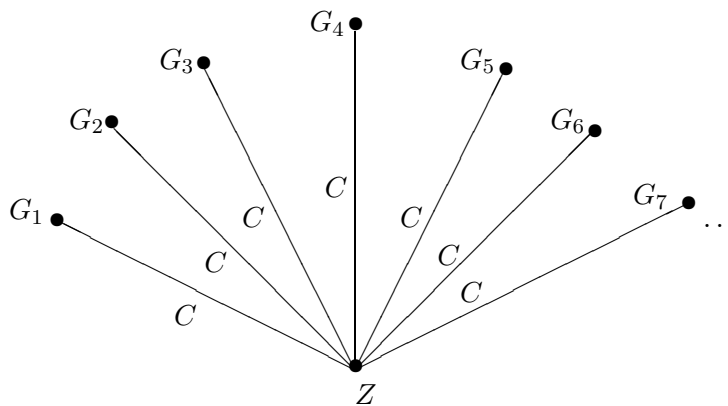
\begin{picture}(100,55)
\put(49,4.6){$\boldsymbol \bullet$}
\put(50,.4){$Z$}
\put(50,5){\line(-2,1){40}}
\put(50,5){\line(-1,1){32}}
\put(50,5){\line(-1,2){20}}
\put(50,5){\line(0,1){45}}
\put(50,5){\line(1,2){20}}
\put(50,5){\line(1,1){32}}
\put(50,5){\line(2,1){45}}
\put(8,24){ $\boldsymbol \bullet$}
\put(26,11){$C$}
\put(4,25){$G_1$}
\put(16.6,37){$\boldsymbol {\bullet}$}
\put(30,18){$C$}
\put(12,37){$G_2$}
\put(28.8,44.8){$\boldsymbol {\bullet}$}
\put(35,24){$C$}
\put(24,45){$G_3$}
\put(49,50){$\boldsymbol {\bullet}$}
\put(46,28){$C$}
\put(44,50){$G_4$}
\put(69,44){$\boldsymbol {\bullet}$}
\put(56,24){$C$}
\put(64,45){$G_5$}
\put(80.8,35.8){$\boldsymbol {\bullet}$}
\put(61,19){$C$}
\put(76,37){$G_6$}
\put(92,26.2){ $\boldsymbol \bullet$}
\put(64,14){$C$}
\put(87,27){$G_7$}
\put(95,24){\dots}
\end{picture}}

\noindent
where each $G_i$ is of the form $G^x$ for some $x \in X$ (where $x$
depends on $i$).  We can now define an
epimorphism $\theta \colon K \twoheadrightarrow G *_C Z$ by
$\theta g = x^{-1}gx$ for $g \in G_i$ and $\theta z = z$ for
$z \in Z$.  The kernel of this map is a free group and hence
left orderable.  Also $G *_C Z$ is left orderable by induction.
We now apply \cite[theorem 7.3.2]{mura}
twice to first deduce that $K$ is left orderable, and then
$G *_C H$ is left orderable, as required.
\end{proof}

\begin{Prob}
Is the free product of two left orderable groups with
an amalgamated cyclic subgroup left orderable?
\end{Prob}

\section{Examples of left ordered groups which are not
locally indicable} \label{Sexamples}

If $G$ is a group, we shall let $\Delta (G)$ indicate the finite
conjugate center of $G$, that is $\{g \in G \mid \Cen_G(g)$ has
finite index in $G\}$.
Let $n$ be a positive integer and let
$B_n$ denote the Braid group on $n$ strings with standard
generators $\sigma_1, \dots, \sigma_{n-1}$.

Dehornoy \cite{dehornoy} (see also \cite{rolfsen})
has proven that the
Braid group $B_n$ on $n$ strings is left orderable.  Therefore $B_n'$
is also left orderable.   It is not difficult to see that
$B_n'$ is a finitely generated
perfect group with trivial center for $n \ge 5$; we shall
give a proof of this (probably known) result in Lemma \ref{Lcenter}
below.  Thus for $n \ge 5$, we see that $B_n'$ is a nontrivial
finitely generated left orderable perfect
group with trivial center; see the last paragraph of
\cite[p.~248]{bergman}.

\begin{Lem} \label{Lfc}
Let $n$ be a positive integer.  Then
$\Delta (B_n) = \Zen(B_n)$.
\end{Lem}
\begin{proof}
Obviously $\Zen(B_n) \subseteq \Delta (B_n)$.  Conversely suppose
$\beta \in \Delta (B_n)$.  Then the centralizer of $\beta$ in
$B_n$ has finite index in $B_n$, consequently it contains a normal
subgroup $C$ of finite index $r$ in $B_n$.
Thus $\sigma_i^r \in C$ for
all $i$, so by \cite[2.2 theorem]{fenn} we see that
$\sigma_i\beta = \beta \sigma_i$ for all $i$.  Therefore $\beta \in
\Zen(B_n)$ and the result is proven.
\end{proof}

\begin{Lem} \label{Lcenter}
Let $n$ be a positive integer.  Then $B_n'$ is
finitely generated and $\Zen(B_n') = 1$.
Furthermore if $n \ge 5$, then $B_n'' = B_n'$.
\end{Lem}
\begin{proof}
The result is trivial if $n \le 2$, so we may assume that
$n \ge 3$.  Let $Z = \Zen(B_n)$.
Then \cite[2.5 corollary]{fenn}
shows that $Z = \langle (\sigma_1 \dots \sigma_{n-1})^n \rangle$,
and we now see from \cite[exercise 7, p.~47]{hansen} that
$Z \cap B_n' = 1$.  Also $B_n/B_n' \cong \mathbb {Z}$ from
\cite[p.~757]{lin} and we deduce that $B_n' Z$ has finite
index in $B_n$.  Therefore $B_n'$ is finitely generated
and $\Zen (B_n') \subseteq \Delta (B_n)$.
But $\Delta (B_n) = Z $ by Lemma~\ref{Lfc} and the first part
is proven.  Finally if $n \ge 5$, then $B_n'' = B_n'$ from
\cite[p.~757]{lin}.
\end{proof}

Let $\tilde {G}$ denote the group of piecewise linear
homeomorphisms of $\mathbb {R}$ which satisfy $g(x+1) = g(x) + 1$
for all $g \in \tilde {G}$ and
$x \in \mathbb {R}$, as described in \cite{ghys}.  Thus
$\tilde {G}$ is a finitely generated
perfect group with infinite cyclic center
$Z$ generated by the map $x \mapsto x+1$ for $x \in \mathbb {R}$,
and $\tilde {G}/Z$ is a simple group, called $T$
in \cite{floyd}.  Then the free product $\tilde {G} *
\tilde {G}$ is a finitely generated
perfect group with trivial center, and is left
orderable by \cite[theorem 7.3.2]{mura}.

\end{document}